\documentclass[12pt,a4paper]{amsart}
\usepackage{amsmath,amssymb,a4}
\usepackage{enumerate,amssymb}
\theoremstyle{plain}
\newtheorem{thm}{Theorem}[section]

\newtheorem{rmk}[thm]{Remark}

\newtheorem{example}[thm]{Example}
\newtoks\prt
\newtheorem{proclaim}[thm]{\the\prt}
\theoremstyle{definition}

\def\eqn#1$$#2$${\begin{equation}\label#1#2\end{equation}}

\numberwithin{equation}{section}

\def\epsilon{\varepsilon}

\def\N{\mathbb N}

\def\loc{\operatorname{loc}}
\def\rn{\mathbb R^n}
\def\spt{\operatorname{supp}}

\def\Ln{\mathcal L^n}
\def\H{\mathcal {H}}
\newtoks\by
\newtoks\paper
\newtoks\book
\newtoks\jour
\newtoks\yr
\newtoks\pages
\newtoks\vol
\newtoks\publ
\def\ota{{\hbox\vol{???}}}
\def\cLear{\by=\ota\paper=\ota\book=\ota\jour=\ota\yr=\ota
\pages=\ota\vol=\ota\publ=\ota}
\def\endpaper{\the\by, {\the\paper},
\the\jour, \the\yr, \the\vol , \the\pages.\cLear}
\def\endbook{\the\by, \the\book, \the\publ.\cLear}
\def\endprep{\the\by, \the\paper, \the\jour.\cLear}
\def\endyearprep{\the\by, \textit{\the\paper}, \the\jour, (\the\yr).\cLear}
\def\name#1#2{#2 #1}

\def\sptt{\operatorname{supp}}
\def\loc{loc}

\def\C{\mathcal{C}}
\def\R{\mathbb {R}}
\def\Rn{\mathbb {R}^n}

\def\applim{\operatorname{app}\ \lim}

\title{Composition operator into the space of function of bounded variation}
\author{Lud\v{e}k Kleprl\'{i}k}
\address{Faculty of Information Technology,
Czech Technical University in Prague, 
Th\' akurova 9, 160 00 Prague 6 
Czech Republic} 
\email{kleprlud@fit.cvut.cz}
\thanks{The author was supported by the grant GA\v CR 18-00960Y.}

\begin{document}
\begin{abstract}
Let $\Omega_1, \Omega_2\subset \R^n$ and $1\leq p <\infty$. We study the optimal conditions on a homeomorphism $f:\Omega_1$ onto $\Omega_2$ which guarantee that the composition $u\circ f$ belongs to the space  $BV(\Omega_1)$ for every $u\in W^{1,p}(\Omega_2)$. We show that the sufficient and necessary condition is an existence of a function $K(y)\in L^{p'}(\Omega_2)$ such that $|Df|(f^{-1}(A))\leq \int_A K(y)\,dy$ for all Borel sets $A$.
\end{abstract}

\maketitle
\section{Introduction}
In this paper we address the following issue. Suppose that $\Omega\subset \Rn$ is an open set, $f : \Omega \to \rn$ is a homeomorphism and a function of bounded variation and  $u$ is a function of $WL^p (f(\Omega))$. Under
which conditions can we then conclude that $u \circ  f \in BV(\Omega)$? Our main theorem gives a complete answer to this question.
\prt{Theorem}
\begin{proclaim}\label{slozen3}
Let $\Omega_1, \Omega_2$ be open subsets of $\R^n$  and let $f\in BV_{\loc}(\Omega_1,\Omega_2)$ have no jump part.   Suppose that $f$ is not constant on any component of $\Omega$ and there is a function $K\in L^{q'}(\Omega_2)$ such that
\eqn{klic55}
$$|Df|(\tilde f^{-1}(A))\leq \int_{A} K(y) \,d \Ln \text{ for all Borel sets } A\subset\Omega_2.$$
Then the operator $T_f(u)(x)=u(f(x))$ maps from function $W^{1,q}(\Omega_2)\cap C(\Omega_2)$ if $q>n$, or   $W^{1,q}(\Omega_2)$ if $q\leq n$ into $BV(\Omega_1)$ and
\eqn{spojitost}
$$|D(u\circ f)|(\Omega_1)\leq \|K\|_{L^{q'}} \|Du\|_{L^q (\Omega_2)}.$$
On the other hand, if $f$ is a homeomorphism of $\Omega_1$ onto $\Omega_2$ such that  the operator $T_f$  maps $C_0(\Omega_2)\cap W^{1,p}(\Omega_2)$  into $BV(\Omega_1)$, 
then $f\in BV_{\loc}(\Omega_1,\Omega_2)$ and there exists a  is a function $K\in L^{q'}(\Omega_2)$ such that \eqref{klic55}  holds.
\end{proclaim}

 The class of homeomorphisms  that satisfy \eqref{klic55} forms a natural extension of a special class of mappings of finite distortion. More precisely: in the fourth chapter we show  that the set of homeomorphisms in $W^{1,1}_{\loc}$ with the property  \eqref{klic55}  coincides with the known class of homeomorphisms with finite distortion satisfying that there exists a function $L\in L^{\frac{1}{q-1}}$ such that
$$|Df(x)|^q\leq L(x)|J_f(x)| \text{ for a.e. } x \in \Omega.$$
It is known that for this class of Sobolev homeomorphisms we have $T_f(u):=u\circ f\in W^{1,1}$ for all $u\in W^{1,q}$. See \cite{GGR} or \cite{VU} for details. Hence naturally $T_f$ maps function from $W^{1,1}$ to $BV$. 

Let us note that the morphism property of $T_f$ on $BV$ was also known under the assumption that the homeomorphism $f$ belongs to class of mappings with a Lipschitz inverse. This can be found in \cite[Theorem 3.16]{AFP}, or \cite{He2}. We show that the above two classes of homeomorphisms differ (and our class contains both of them).

Actually we prove more general statements of the theorems. We allow $f$ to fail to be a homeomorphism. Our mapping will be a general  mapping of bounded variation (its multiplicity can be unbounded) with no jump part and satisfying \eqref{klic55} for some good representative of $f$.

\section{Preliminaries}
We use the usual convention that $C$ denotes a generic positive constant whose exact value may change from line to line.  We denote by $\Ln$ the Lebesgue measure. The symbol $\nabla u(x)$ denotes the classical gradient of $u$ in $x$. By $Du$ we denote the distributional derivative.

Let $\Omega$ be an open subset of $\R^n$. A function $u\in L^1(\Omega)$ whose partial derivatives in the sense
of distributions are measures with finite total variation in $\Omega$ is called a function of bounded variation.  The vector space of functions of bounded variation is denoted by $BV(\Omega)$. We write  $u\in BV(\Omega,\R^d)$ if $u_i \in BV(\Omega)$ for all $i\in\{1,\ldots,d\}$.

If $u\in BV(\Omega,\R^d)$, the total variation of the measure $Du$ is defined by
$$|Du|(E)=\sup\{\sum_{i=1}^m \int_E u_i \operatorname{div} \phi_i \,d\Ln:\phi \in C_c^1(\Omega, \R^{d\times n}), |\phi(x)|\leq 1 \text{ for } x\in \Omega\}<\infty.$$
 We write $u\in BV_{\loc}(\Omega,\rn)$ if for all $x\in \Omega$ there is a ball $B\ni x$ such that $u\in BV(B,\rn)$.
Proposition  3.13 in \cite{AFP} gives us 

\prt{Theorem}
\begin{proclaim}\label{weak}
Suppose that $\Omega\subset \Rn$ is open,  $u_k\in BV(\Omega)$   and there is $u\in L^1(\Omega)$ such that $u_k\rightarrow u$ in $L^1(\Omega)$ and $\sup |Du_k|(\Omega)<\infty$ then $u$ belongs to $BV(\Omega)$ and $u_k$ weakly* converges to $u$ in $BV(\Omega)$.
\end{proclaim}

The main tool is the analogy of the chain rule for the composition of a smooth function and a function of bounded variation see \cite{AM} or Theorem 3.96 in \cite{AFP}.
\prt{Theorem}
\begin{proclaim}\label{dercom}
Suppose that $\Omega\subset \Rn$ is open, $f\in BV (\Omega, \rn)$ and $u\in C^1(\rn,\R^k)$. Then the composition $u\circ f$ belongs to $BV(\Omega)$  and
$$D(u\circ f)= \nabla u\circ f \cdot D^a f \Ln +\nabla u\circ \tilde f\cdot D^c f+ [u(f^+)-u(f^-)]\otimes \nu_f \H^{n-1}|_J,$$
where 
$$Df=D^a f \Ln +D^c f+\nu_f \H^{n-1}|_J$$
is the usual decomposition of $Df$ in its absolutely continuous part $D^a f$ with respect to  the  Lebesgue measure $\Ln$,  its  Cantor  part $D^c u$ and  its  jumping  part,  which  is represented by the restriction of the $(n-1)$ dimensional Hausdorff measure to the jump set $J$.  Moreover, $\nu_f$ denotes the measure theoretical unit normal to $J$, $\tilde f$ is the approximate limit and $f^+$, $f^-$ are the approximate limits from both sides of $J$.
\end{proclaim}
We will work only with functions, which have no jump part, i.e. $J=\emptyset$. In that case we have 
$$D(u\circ f)=\nabla u \circ \tilde f \cdot Df.$$

\subsection{Basic properties of measures}
If $u$ is a  $\mu$-measurable function and $E$ is a $\mu$-measurable set then we denote by $\int_E u\, d\mu$ (or $\int_E u(x)\, d\mu(x)$ if we want to emphasize the variable) the integral of $u$ over $E$ with respect to the measure $\mu$. Instead of $d\Ln(x)$ we write shortly $dx$.

Given measure spaces $(X, \mathcal A)$ and $(Y, \mathcal B)$, a measurable mapping $f : X \to Y$ and a measure $\mu : \mathcal A \to [0, \infty]$, the image of  $\mu$ is defined to be the measure $f(\mu) : \mathcal{B} \to [0, \infty]$ given by
$$(f (\mu)) (A) = \mu \left( f^{-1} (A) \right) \text{ for } A \in \mathcal{B}.$$
Sometimes  $f (\mu)$ is called the pushforward of $\mu$.
\prt{Theorem}
\begin{proclaim}\label{chavar}
Under these assumptions we have that 
\eqn{obrmir}
$$\int_{Y} g \, d(f(\mu)) = \int_{X} g \circ f \, d\mu,$$ 
whenever one of the integrals is well-defined.
\end{proclaim}

Let $\mu, \nu$ be measures defined on the same $\sigma$-algebra $\mathcal A$ of the space $X$. We say that $\mu$  is
\begin{itemize}
\item  absolute continuous with respect to $\nu$ if 
$$|\nu|(A)=0\Rightarrow |\mu|(A)=0.$$
\item singular with respect to $\nu$ if there are $X_a, X_s\in X$ such that $X=X_a\cup X_s$ and
$$|\nu|(X_s)=0=|\mu|(X_a).$$
We set $\sptt \nu=X_s$.
\end{itemize}
For each pair of non-negative measures $\mu$ and $\nu$ on the same $\sigma$-algebra $\mathcal A$ we can find a decomposition $\mu=\mu^a+\mu^c$ such that
$\mu^a$ is absolute continuous with respect to $\nu$ and $\mu^c$, $\nu$ are singular.
\prt{Theorem}
\begin{proclaim}[Radon-Nikodym]\label{radnik}
Let $\mu$ be a non-negative Borel measure on $\rn$ and set 
$$\frac{d\mu}{d\Ln}(x)=\lim_{r\to 0_+} \frac{\mu(B(x,r))}{\Ln(B(x,r))}.$$
Then $\frac{d\mu}{d\Ln}$ exists $\Ln$-a.e., $\frac{d\mu}{d\Ln}(x)$ is $\Ln$-measurable and 
$$\int_A \frac{d\mu}{d\Ln}(x) \, dx \leq \mu(A) \text{ for all Borel sets } A\subset G.$$
Moreover, if $\mu$ is absolute continuous with respect to $\Ln$ then the above inequality holds as equality.
\end{proclaim}

\section{Sufficient condition}
In this section we prove the stability of the composition under our assumptions. We crucially need to know that $f$ satisfies the Lusin $(N^{-1})$ condition, i.e. preimages of sets of Lebesgue measure zero have measure zero. If the condition fails then there is a set $A\subset \Omega_1$ such that $\Ln(A)>0$ and $\Ln\bigl(f(A)\bigr)=0$. 
Now we can redefine $u$ on a null set $f(A)$ arbitrarily and the composed function may be a nonmeasurable function. On the other hand if $f$ satisfies the condition then the validity of our statement for one representative of $u$ implies the validity for all representatives, because the composition only differs on a set of measure zero. The following theorem can be found in \cite{Kl2}
\prt{Theorem}
\begin{proclaim}\label{lusin}
Let $\Omega\subset \R^n$ be connected open set, $f\in BV(\Omega, \rn)$ have no jump part. 
Suppose that
\eqn{key789}
$$|Df|(\tilde f^{-1} (A)) \leq \int_A K(y) \,dy \text{ for all Borel set } A\subset \rn,$$ where $K(y)\in L^{p'}$ for some $p\in [1,n]$, $p'=\frac{p}{p-1}$ . If $f$ is not constant then $f$ satisfies Lusin $(N^{-1})$ condition, i.e. for any set $E\subset\rn$ we have
$$\Ln(E) = 0 \Rightarrow \Ln(f^{-1}(E)) = 0 .$$
\end{proclaim}

\begin{proof}[Proof of the first part of Theorem \ref{slozen3}]
Suppose that $u\in W^{1,q}(\Omega_2)$ . Let be $u_k$ an approximation of $u$ by smooth function in $W^{1,q}(\Omega_2)$ such that $u_k(x)$ converge to $u(x)$ almost everywhere or everywhere if $q>n$. We prove that $u_k\circ f$ is a good weak* approximation of $u\circ f$ on $G$. 
 First by Theorem \ref{dercom} we know that $u_k\circ f$ belongs to $BV(G)$ and $D(u_k\circ f)(x)=\nabla u_k\bigr(f(x)\bigl)\cdot Df(x) $. Hence  with the help Theorem \ref{chavar} and the fact that $\tilde f(|Df|)(A)\leq \int_A K(y) \,  dLn$  we can estimate
\eqn{odh54}
$$
\begin{aligned}
|D(u_k\circ f)|(G)&\leq \int_{G}|\nabla u_k\bigr(\tilde f(x)\bigl)| \,d |Df|=\int_{\Omega_2}|\nabla u_k| \,d (f(|Df|))\\
& \leq  \int_{\Omega_2}|\nabla u_k(y)| K(y) \, d\Ln
\leq  \|K\|_{L^{q'}} \|Du_k\|_{L^q(\Omega_2)}.
\end{aligned}$$ 

 Due to Theorem \ref{lusin} we know  $u_k\circ f\to u\circ f$ almost everywhere on $\Omega_2$ (even everywhere if $q>n$).
Let $B\subset\subset \Omega_1$ be a ball. Due to convergence of $u_k\circ f$ pointwise almost everywhere    using Egorov's theorem we can find a measurable set $A$ such that  $u_k\circ f$ converge to $u\circ f$ uniformly on $A$ and $\Ln(A)=1/2 \Ln(B)$. Denote $v=u_k\circ f-u_l\circ f$  and $v_E=\frac{1}{\Ln(E)}\int_E v(x)\, d\Ln$ then
$$\begin{aligned}
\|v&(x)\|_{L^1(B)}\leq \|v(x)-v_B\|_{L^1(B)}+\Ln(B)(|v_A-v_B|+|v_A|)\\
&\leq \|v(x)-v_B\|_{L^1(B)} + \Ln(B)\left(\left| \frac{1}{\Ln(A)} \int_A v(x)-v_B \, d\Ln\right|+ \left|\frac{1}{\Ln(A)} \int_A v(x)\, d\Ln\right|\right)\\
&\leq \|v(x)-v_B\|_{L^1(B)} +  2 \int_B |v(x)-v_B|\,d\Ln + 2 \|v\|_{L^\infty(A)}\\
&\leq  3 r|Dv|(B) + 2 \|v\|_{L^\infty(A)},
\end{aligned}$$
where we used the Poincare inequality. Thus together with \eqref{odh54} we have
$$\|u_k\circ f-u_l\circ f\|_{L^1(B)}\leq 3  \|K\|_{L^{q'}} \|Du_k-Du_l\|_{L^q(\Omega_2)} + 2 \|u_k\circ f-u_l\circ f\|_{L^\infty(A)}.$$
It follows $u_k\circ f$ is a  Cauchy sequence in $L^1(B)$ for all balls. Due to almost everywhere pointwise convergence of $u_k\circ f$ we have that $u_k\circ f$ converge to $u_k\circ f$ in $L^1_{\loc}(\Omega_1)$.

Thus convergence of $u_k\circ f$ to $u\circ f$ in $L^1(G)$ together with the estimate \eqref{odh54} and Lemma \ref{weak} give us that  $u\circ f$ is function of bounded variation on $G$.

Moreover,   using semi-continuity of the variation we obtain
$$
\begin{aligned}
|D(u\circ f)|(G)&\leq \liminf_{k\to \infty} |D(u_k\circ f)|(G)\leq \liminf_{k\to \infty} \|K\|_{L^{q'}} \|u_k\|_{L^q(\Omega_2)}=  \|K\|_{L^{q'}} \|u\|_{L^q(\Omega_2)}.
\end{aligned}$$
To prove \eqref{spojitost} find open sets $G_k\subset\subset \Omega$ such that $G_k\subset G_{k+1}$ and $\Omega_1=\bigcup_{1}^{\infty} G_k$ then 
$$|D(u\circ f)|(\Omega_1)=\lim_{k} |D(u\circ f)|(G_k)\leq \|K\|_{L^{q'}} \|u\|_{L^q(\Omega_2)}.$$

\end{proof}
In the case when $f$ is constant on some component $G$ of $\Omega$ the composition $u\circ f$ may not be well-defined. If we take a representative of $u$ such that $\tilde u(x)=0$ for all $x$ such there is a component $G$ of $\Omega$ satisfying $f(G)=\{x\}$ then for this representative we have $\tilde u\circ f \in BV(\Omega_1)$ and  \eqref{spojitost} again holds.

\begin{rmk}
The condition \eqref{klic55} can be rewritten
\eqn{bobabobek}
$$\int_{\tilde f^{-1} (A)} |D^a f|\,d\Ln + \int_{\tilde f^{-1} (A)} \,d|D^c f| \leq \int_{A} K(y) \,d\Ln ,$$
what is equivalent to a existence of function  $K^a,K^c \in L^q$ such that 
\eqn{blb}
$$\int_{\tilde f^{-1} (A)} |D^a f|\,d\Ln \leq \int_{A} K^a(y) \,d\Ln$$
 and 
\eqn{bob}
$$ \int_{\tilde f^{-1} (A)} \,d|D^c f|\leq \int_{A} K^c(y) d\,\Ln.$$
The second condition \eqref{bob} implies that $|D^c f|(\tilde f^{-1}(A))=0$ whenever $A\subset \Omega_2$ has measure zero.
\end{rmk}

\prt{Lemma}
\begin{proclaim}
Assume that $f$ is a homeomorphism then the inequality \eqref{blb} is equivalent to 
existence of a function $L\in L^{q'}$ such that 
\eqn{1wconf}
$$|D^a f(x)| \leq L(x) |J_f(x)|^{\frac{1}{q}} \text{ for a.e. }x \in \Omega_1.$$
\end{proclaim}
\begin{proof}
 Find Borel  set $Z$ of measure zero such that  $f|_{\Omega_1\setminus Z}$ satisfies the Lusin $(N)$ condition and set $N=\{J_f=0\}$.

 If \eqref{1wconf} holds then using the Area formula we obtain 
$$\begin{aligned}
\int_{ f^{-1} (A)} |D^a f|\,d\Ln&=\int_{ f^{-1} (A)\setminus (N\cup Z)} |D^a f|\,d\Ln\leq \int_{ f^{-1} (A)\setminus (N\cup Z)} L(x) |J_f(x)|^{\frac{1}{q}} \,d\Ln\\
&= \int_{A\setminus f(N\cup Z)} L(f^{-1}(y)) |J_f(f^{-1}(y))|^{\frac{1-q}{q}}\,d\Ln
\end{aligned}$$
So we can set $K(y)= L(f^{-1}(y)) |J_f(f^{-1}(y))|^{\frac{1-q}{q}} \chi_{\rn\setminus f(N\cup Z)}$
where
 $$\int_{A} \left(L(f^{-1}(y))^{q'} |J_f(f^{-1}(y))|^{\frac{1-q}{q}} \chi_{\rn\setminus f(N\cup Z)}\right)^{q'}\,d\Ln= \int_{f^{-1}(A)} L(x)^{q'} \,d\Ln$$

Let us assume that \eqref{bobabobek} hold. Then  
 $$\begin{aligned}
\int_{ B(x,r)} |D^a f|\,d\Ln&=\int_{ B(x,r)\setminus Z} |D^a f|\,d\Ln\leq  \int_{f(B(x,r)\setminus Z)} K(y) \,d\Ln\\
&\leq \left(\int_{f(B(x,r)\setminus Z)} K(y)^{q'} \,d\Ln\right)^{\frac{1}{q'}} \left(f(B(x,r)\setminus Z)\right)^{\frac{1}{q}}\\
&= \left(\int_{(B(x,r)\setminus Z)} K(f(x))^{q'} |J_f(x)| \,d\Ln\right)^{\frac{1}{q'}} \left(\int_{B(x,r)} |J_f(x)|\right)^{\frac{1}{q}}.
\end{aligned}$$
Considering that $K(f(x))^{q'} |J_f(x)|$ is integrable then by letting $r\to 0_+$ we obtain for all Lebesgue points
$$|D^a f(x)|\leq \left(K(f(x))|J_f(x)|^{\frac{1}{q'}}\right) |J_f(x)|^{\frac{1}{q}}.$$  It can be easily verified that $L(x)=K(f(x))|J_f(x)|^{\frac{1}{q'}}\in L^{q'}$.
 \end{proof}

If we assume that $f$ is a Sobolev homeomorphism then $D^c f=0$. 
\prt{Corollary}
\begin{proclaim}
Let $f$ be a homeomorphism in $W_{\loc} ^{1,1}(\Omega_1,\R^n)$ then \eqref{klic55} is equivalent to
existence of a function $L\in L^{q'}(\Omega_1)$ such that 
\eqn{1conf}
$$Df(x)\leq L(x) |J_f(x)|^{\frac{1}{q}} \text{ for a.e. }x \in \Omega_1.$$
\end{proclaim}

The simplest way to obtain the condition \eqref{klic55} is to check the integrability of the inverse.
\prt{Lemma}\begin{proclaim}\label{lip2}
Let $\Omega_1,\Omega_2 \subset \rn$ and $f:\Omega_1\to \Omega_2$ be a mapping such that $f^{-1}$ is in $L^p$, where $p$ is bigger or equal to $n-1$. Then \eqref{klic55} holds for $q=\frac{p}{p-(n-1)}$.
\end{proclaim}
\begin{proof}
It follows from Lemma 4.3 in \cite{CHM} and Theorem 3.8 \cite{Mi} that  $f\in BV_{\loc}(\Omega_2)$ and 
\eqn{odh11}
$$|Df|(f^{-1}(G))\leq C\int_{G} |\operatorname{adj} D(f^{-1})| \,d\Ln,$$ where $C$ depends only on $n$. Hence \eqref{odh11} holds for all Borel sets and we have 
$$|Df|(f^{-1}(A))\leq  C\int_{A} |\operatorname{adj} D(f^{-1})| \, d\Ln.$$ Taking in mind $p=\frac{q(n-1)}{q-1}$ we can estimate
$$  \int_{\rn} \operatorname{adj} |D(f^{-1})^{q'}| \, d\Ln\leq \int_{\rn} \operatorname{adj} |D(f^{-1})|^p \, d\Ln.$$
\end{proof}
\begin{example}
There is a homeomorphism $f$ such that \eqref{klic55} holds even for $q=1$ but $f^{-1}\notin W^{1,1}_{\loc}$.
\end{example}
\begin{proof}
 Consider the usual Cantor ternary function $u$ on the interval $(0,1)$. And set $g(x)=u(x)+x$. This function is continuous, increasing and fails to be absolutely continuous. Moreover, $g$ does not belong to $W^{1,1}_{\loc}$. On the other hand, the inverse function $g^{-1}$ is Lipschitz and maps $(0,2)$ homeomorphically onto $(0,1)$. If we set 
$$f(x_1,\ldots,x_n)=(g(x_1),x_2,\ldots, x_n)$$ then obviously $f$ fails to belong to $W^{1,1}_{\loc}((0,1)^n)$, and $f^{-1}$ is a  Lipschitz function. Due to Lemma \ref{lip2} the function $f$ satisfy \eqref{klic55}.
\end{proof}

In the special case when $n=2$ we obtain the equivalence in Lemma \ref{lip2}. 
\prt{Lemma}\begin{proclaim}
Let $\Omega_1,\Omega_2 \subset \R^2$ and $f:\Omega_1\to \Omega_2$ be a homeomorphism. Then $f^{-1}\in W^{1,q}(\Omega_2)$ if and only if $f\in BV_{\loc}(\Omega_1)$ and \eqref{klic55} holds.
\end{proclaim}
\begin{proof}
It remains to prove the second implication. Let $f\in BV_{\loc}(\Omega_1)$. It follows from \cite{DS} 
that $f^{-1}$ is in $BV_{\loc}(\Omega_2)$ and
$$|D(f^{-1}_1)|(\Omega_2)=|D_y f|(f^{-1}(\Omega_2)) \text{ and }|D(f^{-1}_2)|(\Omega_2)=|D_x f|(f^{-1}(\Omega_2)).$$
It holds for all open set $\Omega_2$ thus we have for all Borel sets $A\subset \Omega_2$
$$|D(f^{-1})|(A)= |Df|(f^{-1}(A)).$$
By combining with \eqref{klic55} we have $|D(f^{-1})|(A)\leq \int_A K(y) \,d \Ln$ and thus the measure $D(f^{-1})$ is absolutely continuous with respect to $\Ln$ and which Radon-Nikodym derivative with respect to $\Ln$ is in $L^{q}$.
\end{proof}

\section{A necessary condition}
In this section we prove the second part of the main Theorem \ref{slozen3}
\prt{Theorem}
\begin{proclaim}\label{neces}
Let $\Omega_1, \Omega_2$ be open subsets of $\R^n$  and let $f\in BV_{\loc}(\Omega_1,\Omega_2)$ has no jump part and suppose that  the operator $T_f$ maps functions from $C_c^\infty(\Omega_2)$ into $BV_{\loc}(\Omega_1)$ and there is a $M\in \R$ such that for all $u\in C_c^\infty(\Omega_2)$ we have
\eqn{necesklic}
$$|D(u\circ f)|(\Omega_1)\leq M \|Du\|_{L^q (\Omega_2)}.$$
 Then  there is a function $L(y)\in L^{\frac{q}{q-1}}$ such that $\int L(y)^{\frac{q}{q-1}} \,d \Ln\leq 16n M$ for all Borel set $A\subset \Omega_2$  we have
\eqn{klic56}
$$|Df|(\tilde{f}^{-1}(A))\leq \int_A {L(y)} \,d \Ln.$$
\end{proclaim}
\begin{proof}
For  an open set $G\subset \Omega_2$ we denote by $M(G)$ the smallest constant such that \eqref{necesklic} holds for all $u\in C_c^\infty(G)$. Then obviously $M$ is a set function on open sets. First note that $M(G)$ (and  thus also $M(G_k)^{\frac{q}{q-1}}$) is monotone i.e. $M(G)\leq M(\tilde G)$ if $G\subset \tilde G$. We claim that $M(G_k))^{\frac{q}{q-1}}$ is quasiaditive i.e. for any finite collection $G_1,\ldots G_k$ of pairwise disjoint open subsets of $\Omega_2$ we have  
$$M(\bigcup_k G_k)^{\frac{q}{q-1}}\leq \sum M(G_k)^{\frac{q}{q-1}}.$$
To see it take $u\in C^\infty_c(\bigcup G_k)$. Then $u_k=u \chi_{G_k} \in C_c^\infty(G_k)$. From Theorem \ref{dercom} we can see that $|Du_k\circ f|(\Omega_1\setminus f^{-1}(G_k))=0$ thus $|D(u\circ f)|=\sum |D(u_k\circ f)|$ and
$$\begin{aligned}
|D&(u\circ f)|(\bigcup_k G_k)\leq \sum_k M(G_k)\|Du_k\|_{L^q(G_k)}
\\&\leq \left(\sum_k (M(G_k))^{\frac{q}{q-1}}\right)^{\frac{q-1}{q}} 
\left(\sum_k \|Du_k\|^q_{L^q(G_k)}\right)^{\frac{1}{q}}\\
&= \left(\sum_k (M(G_k))^{\frac{q}{q-1}}\right)^{\frac{q-1}{q}} \|Du_k\|_{L^q(\bigcup_k G_k)}.
\end{aligned}$$
Thus $$M(\bigcup_k G_k)\leq \left(\sum M(G_k)^{\frac{q}{q-1}}\right)^{\frac{q-1}{q}}.$$


We may assume that $f=\tilde f$. (We change $u\circ f$ only on a set of measure zero.) Take $G\subset \Omega_2$ an open set. 
  First suppose that  $|D f|(f^{-1}(G))\neq 0$. 
	Let $t>0$ and $0<L< |Df|(f^{-1}(A))$ be  arbitrary real numbers and fix $i\in \{1,\ldots,n\}$ such that $|D(f_i)|(A)\geq \frac{1}{n} |Df|(f^{-1}(A))>\frac{1}{n} L $. 
 Then $$G=\bigcup_k A_k=\bigcup_k \{x\in G\cap B(0,k): \operatorname{dist}(x,\partial G)\geq 1/k\}.$$
Choose $k\in \N$ big enough such that $$|Df|(f^{-1}(A_k))>\frac{1}{n} L.$$
Find  a cut-off function $\eta\in \C_c^{\infty}(\Omega_2)$ 
satisfying  
$$\spt\eta\subset G,\ 0\le \eta\le 1\text{ and }\eta=1 \text{ on }A_k.$$ 
Take $m$ such that $m\geq 8$ and
$\|\nabla \eta\|_{\infty}\leq m$. 
Choose $K$ among the sets
$$
\aligned
K^{\sin}&=
\{x\in f^{-1}(A_k):\ \cos^2 ( m^2  f_i(x))\geq \tfrac12\},\\
K^{\cos}&=
\{x\in f^{-1}(A_k):\ \sin^2 ( m^2  f_i(x))\geq \tfrac12\}
\endaligned
$$
such that
$$
|D(f_i)|(K)\geq \tfrac 12 |D(f_i)|(f^{-1}(A_k))
$$
and set 
$$
u (y)=
\begin{cases}
\frac{1}{m^2} \eta(y)\sin (m^2 y_i)&\text{ if }K=K^{\sin}
\\
\frac {1}{m^2} \eta(y)\cos (m^2 y_i)&\text{ if }K=K^{\cos}.
\end{cases}
$$
Consider $K=K^{\sin}$. The case  when $K=K^{\cos}$ is analogous. Obviously $u\in C^\infty_c(\Omega_2)$ and 

\eqn{lips}
$$|\nabla u(y)|= |1/m^2\nabla \eta(y) \sin (m^2 y_i)+\eta(y) \cos(m^2 y_i) e_i|\leq 2\text{ for all } y\in \Omega_2.$$

By the  product rule  it easily follows 
$$
\begin{aligned}
|D(u\circ f)|(K)&= \int_K  \,d|D(u \circ f)|=\int_K  |\nabla u|(\tilde f(x)) \,d|D f_i|\\
&\geq \int_K \left(|\eta(y) \cos(m^2 y_i) e_i|-|1/m^2\nabla \eta(y) \sin (m^2 y_i)|\right)\, d|D f_i|\\
&\geq \int_K (\tfrac{1}{\sqrt{2}}-\tfrac{1}{m})\, d|D f_i| \geq \frac 1{4}\ |D f_i|(K)\\
&\geq \frac 18|Df_i|(f^{-1}(A_k))\geq \frac{1}{8n} L.
\end{aligned}$$

Thus together with \eqref{necesklic} we estimate
$$\begin{aligned}
 L\leq  8n |D(u\circ f)|(\Omega_1)&\leq  8n\, M(G) \|Du\|_{L^q (\Omega_2)} \leq 8n  M(G)\, 2\cdot(\Ln(G))^\frac{1}{q}\leq 16n M(G)  (\Ln(G))^\frac{1}{q}.
\end{aligned}$$
By taking supremum over all $L\leq |D f|(f^{-1}(A))$  we obtain 
\eqn{mg}
$$|D(u\circ f)|(f^{-1}(G))\leq 16n M(G)  \Ln(G)^\frac{1}{q}.$$
Obviously \eqref{mg} is valid even if $|D(u\circ f)|(f^{-1}(G))=0$. Set $\mu(A):=f(|Df|)(A)=|D(u\circ f)|(f^{-1}(A))$ if $A$ is Borel. Because $M(G)\leq M$ we can estimate 
\eqn{measure}
$$\mu(G)\leq 16n M  (\Ln(G))^\frac{1}{q}.$$
Thus $\mu$ is absolute continuous with respect to the Lebesgue measure. Now we apply \eqref{mg} on $B=B(y,r)$ to get 
$$\frac{\mu(B)}{\Ln(B)}\leq 16n \left(\frac{M(B)^{\frac{q}{q-1}}}{\Ln(B)}\right)^{\frac{q-1}{q}}.$$
We have proved in that $\Psi(G)=(M(G))^{\frac{q}{q-1}}$ is a monotone and quasiaditive set function of open sets. Thus 
$$L(y):=\frac{\partial \mu}{\partial \Ln}(y)\leq  16n \left(\Psi'(y)\right)^{\frac{q-1}{q}} \text{ for a.e. }y.$$
Integrating over $G$ taking into account absolute continuity of $\mu$ we obtain
\eqn{konec1}
$$|D(u\circ f)|(f^{-1}(G))=\mu(G)= \int_G L(y)$$
and  
$$\int_G L(y)^{\frac{q}{q-1}}\leq \int_G 16n \Psi'(y) \leq 16n \Psi(G)=16n M.$$
\end{proof}

\prt{Lemma}
\begin{proclaim}\label{pulim}
Let $G$ be an open subset of $\rn$ and $\mu:\mathcal B(G)\to [0,\infty]$ be a non-atomic Borel measure. Then there are open pairwise disjoint sets $G_1,G_2\subset G$ such that
$\mu(G_1)=\mu(G_2)$.
\end{proclaim}
\prt{Lemma}
\begin{proclaim}\label{unb}
Let $G$ be an open subset of $\rn$ and $g$ be a $\Ln$-measurable and essentially unbounded function. Then there are open pairwise disjoint sets $G_k\subset G$, $k\in \N$ such that $g|_{G_k}$ are essentially unbounded.
\end{proclaim}

\begin{proof}[Proof of the second part of Theorem \ref{slozen3}]
Suppose that \eqref{klic55} is not satisfied. Then there are two cases

1. The measure $\nu=f(|Df|)$ is not absolute continuous with respect to $\Ln$. Denote by $\mu$ the singular part. Then with the help of Lemma \ref{pulim} we can  find pairwise disjoint open sets $G_k$ such that $\mu(G_k)>0$.

2.  The measure $\nu$ is absolute continuous with respect to $\Ln$. It means that $\frac{\partial \nu}{\partial \Ln}\notin L^{p'}$. If $p>1$ consider measure $\mu(A)= \int_A \left(\frac{\partial \nu}{\partial \Ln}\right)^{p'}\, d\Ln$. Using Lemma \ref{pulim} we find pairwise disjoint open sets $G_k$ such that $\mu(G_k)=\infty$.

If $p=1$ the Lemma \ref{unb} gives us pairwise disjoint open sets $G_k$ such that $\frac{\partial \nu}{\partial \Ln}$ is essentially unbounded on $G_k$.

In both cases we have that \eqref{klic55}  is not satisfied on $G_k$ hence the assumptions
of Theorem $\ref{neces}$ cannot be satisfied there. It follows that we can construct $u_k\in C^\infty_C (G_k))$ such that
$$\|Du_k\|_{L^q(G_k)}\leq 2^{-k}$$ and
$$|D(u_k\circ f)|(f^{-1}(G_k))\geq 2^k. $$ 
\end{proof}
Due to the fact that $D|v|=|Dv|$ for any function $v$ of bounded variation or Sobolev function we may assume that $\|u_k\|_{L^\infty} \leq 1$ (Otherwise we can  iterate replacing $u_k$ by function $\tilde{u_k}=||u_k|-1/2\|u_k\|_{L^\infty}|$, which has the same total variation of the distributional derivative and its maximum is half of the maximum of $u_k$.)
We extend the domain of the functions $u_k$ by putting $u_k = 0$ on $\Omega\setminus G_k$. Set
$$u = \sum^\infty_{k=1} u_k .$$
Then the sums converges in the norm of $W_0^{1,q} (\Omega)$ and $C_0$. Now, assume that the function
$u \circ f$ is a Sobolev function on $\Omega$, otherwise there is nothing to prove. Then $u_k\circ f$ has a compact support in $f^{-1}(G_k)$ thus 
$$|D(u \circ f )|(\Omega) \geq |D(u \circ f )|(f^{-1}(G_k))= |D(u_k \circ f )|(f^{-1}(G_k))\geq 2^k.$$
It easily follows that $u \circ f  \notin BV(\Omega)$.

\end{document}